\documentclass{article}
\PassOptionsToPackage{hidelinks}{hyperref}
\usepackage[preprint]{log_2026}			

\usepackage{booktabs}						
\usepackage{multirow}						
\usepackage{amsfonts}						
\usepackage{graphicx}						
\usepackage{duckuments}						

\usepackage[numbers,compress,sort]{natbib}	

\usepackage{amsthm}
\usepackage{thm-restate}

\usepackage{amsmath}
\usepackage{subcaption}
\usepackage{multirow}
\usepackage{tikz}
\usepackage{tikz-cd}
\usetikzlibrary{calc,positioning}
\newcommand{\R}{\mathbb{R}}
\newcommand{\Z}{\mathbb{Z}}
\newcommand{\N}{\mathbb{N}}
\AtEndPreamble{\RequirePackage[nameinlink]{cleveref}}
\usepackage{comment}

\newtheorem{theorem}{Theorem}[section]

\newtheorem{definition}[theorem]{Definition}

\newtheorem{lemma}[theorem]{Lemma}
\newtheorem{remark}[theorem]{Remark}

\usepackage[utf8]{inputenc} 
\usepackage[T1]{fontenc}    
\usepackage{url}            
\usepackage{booktabs}       
\usepackage{amsfonts}       
\usepackage{xcolor}         

\title[Stable Vectorization of Persistent Laplacians via Spectral Descriptors]{Stable Vectorization of Persistent Laplacians via Spectral Descriptors}

\author[Jung et al.]{%
  Inkee Jung\thanks{Equal contribution.}\\
  Department of Mathematics and Statistics\\
  Boston University\\
  Boston, MA, USA \\
  \texttt{inkeej@bu.edu}\And
  Wonwoo Kang\footnotemark[1]\\
  International Center for Mathematical Sciences \\ Institute of Mathematics and Informatics\\
  Bulgarian Academy of Sciences\\
  Sofia, Bulgaria \\
  \texttt{wonk@math.bas.bg}\And
  Heehyun Park\footnotemark[1]\\
  Department of Computer Science and Engineering \\ Pennsylvania State University \\
  State College, PA, USA \\
  \texttt{hbp5148@psu.edu}
}

\begin{document}

\maketitle

\begin{abstract}
Persistence images vectorize persistence diagrams into stable, finite-dimensional features. Inspired by this idea, we developed a vectorization framework for the spectral information encoded by the persistent Laplacian (PL). Given a scalar \emph{signature} of a persistent Laplacian, we form a \emph{Persistent Laplacian Diagram} (PLD) and smooth it into a \emph{Persistent Laplacian Image} (PLI). We prove a stability theorem for PLIs with respect to the Wasserstein distance between PLDs under an admissibility condition on the signature. Through experiments on MNIST and QM7, we show that PLI with suitable signatures, especially the trace, provides an effective way to extract predictive topological and geometric information from PL, outperforming existing PL-based representations in these settings.
\end{abstract}

\section{Introduction}
\label{sec:intro}
\emph{Persistent Homology} (PH) plays a central role in topological data analysis (TDA). It tracks how homological features appear and merge along a filtration, producing multi-scale descriptor \emph{Persistence Diagrams} (PD). PDs are multisets of points in the birth--death plane that are stable under perturbations of the filtration \cite{CEH07, ELZ02}. However, at the same time, because they are not fixed-dimensional vectors, PDs often need to be transformed into representations better suited for machine learning. \emph{Persistence Images} (PI)~\cite{PI17} are standard examples: they transform persistence diagrams into stable finite-dimensional features by smoothing and discretization~\cite{PI17}.
 
Beyond homology groups, \emph{Persistent Laplacian} (PL) associates a combinatorial Laplacian with pairs of steps in a filtration~\cite{davies2023persistent,memoli2022persistent,WNW20,WW25}. By construction, the kernel of the persistent Laplacian recovers the persistent homology while containing the nonzero spectrum, which captures homotopic shape evolution and geometric information not detected by homology \cite{chen2022persistent, memoli2022persistent}. Persistent Laplacians have been used in many applications, including biomolecule data and analyzing point cloud data \cite{chen2022persistent,wei2023persistent,WW25}. However, despite this richer information, there is still no broadly used PL analogue of persistence images that is simultaneously stable and fixed-dimensional. This limits the use of persistent Laplacians as generic features for downstream learning tasks.


We address this gap by introducing a PI-style vectorization framework for persistent Laplacians. Given a scalar signature of a persistent Laplacian, we first form a persistent Laplacian diagram (PLD), which records how the chosen signature varies across birth--death pairs. We then smooth and discretize this diagram to obtain a persistent Laplacian image (PLI), a fixed-dimensional representation that can be used directly in standard learning pipelines. This construction is modular in the choice of signature and parallels the role of persistence images for persistence diagrams. We also identify conditions on the signature under which the resulting PLI is stable with respect to Wasserstein perturbations of the PLD.

Experiments on the MNIST and QM7 datasets demonstrate that the proposed PLI framework outperforms the reproduced PL baselines established in \cite{davies2023persistent}. Specifically, we observe that the choice of spectral signature acts as a crucial parameter that affects performance. Among the signatures, the trace yields superior performance across both datasets.


This paper makes four main contributions. We introduce a signature-based vectorization framework for persistent Laplacians through persistent Laplacian diagrams and persistent Laplacian images. We show how signatures can be incorporated into this framework, yielding flexible representations with expressive properties. We establish stability for persistent Laplacian images with respect to the Wasserstein distance between PLDs. Finally, we evaluate the proposed representations on MNIST and QM7, where the trace-based PLI achieves the strongest performance among the signatures considered.




\section{Preliminaries}
\label{sec:prelim}

We introduce preliminaries for our paper: persistent homology, persistent Laplacian, and persistence diagram. We adopt the following notations and definitions from~\cite{memoli2022persistent, WW25}.

\subsection{Persistent Homology(PH) and Persistent Laplacian(PL)}

Both persistent homology and persistent Laplacians are built upon the fundamental structure of chain complexes. In this subsection, we introduce the notion of a chain complex and describe how persistent homology and persistent Laplacians are constructed from it.

A \emph{simplicial complex} $K$ over a finite ordered set $V$ is a collection of finite subsets of $V$ satisfying the condition that for any $\sigma \in K$, if $\tau \subseteq \sigma$, then $\tau \in K$. Let $\mathbb{N}$ denote the set of non-negative integers. For each $q \in \mathbb{N}$, an element $\sigma \in K$ is called a \emph{$q$-simplex} if $|\sigma| = q + 1$, where $|\sigma|$ denotes the cardinality of $\sigma$. A $0$-simplex, typically denoted by $v$, is referred to as a \emph{vertex}. Denote by $S_q^K$ the set of $q$-simplices of $K$, noting that $S_0^K \subseteq V$. The \emph{dimension} of $K$, denoted $\dim(K)$, is the largest integer $q$ such that $S_q^K \neq \emptyset$. A simplicial complex of dimension $1$ is also called a \emph{graph}, and we often represent it as $K = (V^K, E^K)$, where $V^K = S_0^K$ is the vertex set and $E^K = S_1^K$ is the edge set.

An \emph{oriented simplex}, denoted by $[\sigma]$, is a simplex $\sigma \in K$ equipped with an ordering of its vertices. We always assume that this ordering is inherited from the ordering of $V$. Let $\bar{S}_q^K := \{[\sigma] : \sigma \in S_q^K\}$. The $q$-th \emph{chain group} of $K$, denoted $C_q^K := C_q(K, \mathbb{R})$, is the vector space over $\mathbb{R}$ spanned by $\bar{S}_q^K$. Its dimension is given by $n_q^K := \dim C_q^K = |S_q^K|$. The \emph{boundary operator} $\partial_q^K : C_q^K \to C_{q-1}^K$ is defined by

\(
\partial_q^K([v_0, \ldots, v_q]) := \sum_{i=0}^q (-1)^i [v_0, \ldots, v_{i-1}, v_{i+1}, \ldots, v_q]
\)
for each oriented simplex $[\sigma] = [v_0, \ldots, v_q] \in \bar{S}_q^K$.
The \emph{$q$-th homology group} of $K$ is defined as \(H_q(K) = \frac{\ker(\partial_q^K)}{\text{Im}(\partial_{q+1}^K)},\) and its rank $\beta_q^K := \text{rank}(H_q(K))$ is called the \emph{$q$-th Betti number}.

A \emph{weight function} on a simplicial complex $K$ is any positive function $w^K : K \to (0, \infty)$. Throughout this paper, every simplicial complex $K$ is implicitly assumed to be endowed with such a weight function.
Let $K$ be a simplicial complex with a weight function $w^K$. Given any $q \in \N$, let $w^K_q:=w^K|_{S^K_q}$. We define an inner product $\langle\cdot,\cdot\rangle$ on $C_q^K$ by
\[
\langle [\sigma], [\sigma'] \rangle_{w_q^K} := \delta_{\sigma \sigma'} \left(w_q^K(\sigma)\right)^{-1}, \quad \forall\, \sigma, \sigma' \in S_q^K.
\]

Let $\left(\partial_q^K\right)^* : C_{q-1}^K \to C_q^K$ denote the adjoint of $\partial_q^K$ with respect to this inner product. The \emph{$q$-th combinatorial Laplacian} $\Delta_q^K : C_q^K \to C_q^K$ is then defined as
\[
\Delta_q^K := \partial_{q+1}^K \circ \left(\partial_{q+1}^K\right)^* + \left(\partial_q^K\right)^* \circ \partial_q^K.
\]
The first term of $\Delta_q^K$ is referred to as the \emph{up Laplacian}, and the second as the \emph{down Laplacian}. Setting $\partial_0^K = 0$ yields $\Delta_0^K = \partial_1^K \circ \left(\partial_1^K\right)^*$. When $K$ is a graph and $w_0^K = 1$, $\Delta_0^K$ coincides with the classical graph Laplacian of the weighted graph $(K, w_1^K)$~\cite{Chung97}.

A \emph{simplicial pair}, denoted $K \hookrightarrow L$, consists of two simplicial complexes $K$ and $L$ defined on the same finite ordered vertex set $V$, satisfying $K \subseteq L$—that is, $S_q^K \subseteq S_q^L$ for all $q \in \mathbb{N}$ and $w^K = w^L|_K$. A \emph{simplicial filtration} $\mathbf{K} = \{K_t\}_{t \in T}$ is a family of simplicial complexes over a common vertex set $V$, indexed by a subset $T \subseteq \mathbb{R}$, such that for all $s \leq t$ in $T$, the inclusion $K_s \hookrightarrow K_t$ forms a simplicial pair. Another way of understanding subscripts in $K_t$ is, we can set a real function $f:K\rightarrow \R$ where $K=\cup_{t\in T}K_t$ and $K_t=f^{-1}((-\infty,t])$.

For any integer $q \ge 0$ and for all $s \le t \in T$, the inclusion $K_s \hookrightarrow K_t$ induces a homomorphism \(f_q^{s,t} : H_q(K_s) \to H_q(K_t)\). Now, suppose we have a simplicial pair $K \hookrightarrow L$ and fix $q \in \mathbb{N}$. Define the subspace
\[
C_q^{L,K} := \{\, c \in C_q^L : \partial_q^L(c) \in C_{q-1}^K \,\} \subseteq C_q^L,
\]
which consists of $q$-chains in $C_q^L$ whose boundaries lie in the subspace $C_{q-1}^K \subseteq C_{q-1}^L$.  
For each $q$, let $\partial_q^{L,K}$ denote the restriction of $\partial_q^L$ to $C_q^{L,K}$, giving rise to the ``diagonal’’ boundary operator \(\partial_q^{L,K} : C_q^{L,K} \to C_{q-1}^K\). See the diagram below for the construction.

\begin{center}
    
\begin{tikzcd}[scale=0.8, transform shape]
  C_{q+1}^K \arrow[hookrightarrow, gray]{dd} \arrow{rr}{\partial_{q+1}^K} & & C_{q}^K \arrow[hookrightarrow, gray]{dd} \arrow[blue]{rr}{\partial_{q}^K} &   & C_{q-1}^K \arrow[hookrightarrow, gray]{dd} \\
  & C_{q+1}^{L,K}\arrow[hookrightarrow, gray]{dl} \arrow[blue]{ur}{\partial_{q+1}^{L,K}} & & & \\
  C_{q+1}^L\arrow{rr}{\partial_{q+1}^L} & &  C_q^L \arrow{rr}{\partial_{q}^L}&  & C_{q-1}^L
\end{tikzcd}
\end{center}

The adjoint $\left(\partial_q^{L,K}\right)^*$ of $\partial_q^{L,K}$, is taken with respect to the inner product $\langle \cdot, \cdot \rangle_{w_q^L}$. We define the \emph{$q$-th persistent Laplacian}~\cite{WNW20} $\Delta_q^{L,K} : C_q^K \to C_q^K$ by
\[
\Delta_q^{L,K} := \partial_{q+1}^{L,K} \circ \left(\partial_{q+1}^{L,K}\right)^* 
+ \left(\partial_q^K\right)^* \circ \partial_q^K.
\]
It is immediate that when $K = L$, the operator $\Delta_q^{L,K}$ reduces to $\Delta_q^L$, the standard $q$-th Laplacian on $L$.

\subsection{Persistence Diagrams and Persistence Image}

A standard way to represent PH information is a persistence diagram~\cite{ELZ02} or barcodes~\cite{Ghrist07}. Here, we introduce \emph{persistence diagram}. Some of these points may have infinite coordinates or coincide with others. Thus, we consider a multiset of points in the extended real plane $\overline{\R}^2$, where $\overline{\R} = \R \cup \{\infty\}$. Since each point has the form $(b,d)$ with $d>b$, we take $\mathcal{R} = \R \times \overline{\R}$ as the ambient set of persistent diagrams. For $p = (b,d) \in \{(u,v) \in \mathcal{R} : u < v\}$, let $\mu_{q}^{b,d}$ denote the number of $q$-dimensional homology classes that are born at $K_b$ and die upon entering $K_d$. Then, for all $b < d$ and all $q \in \Z_{\geq 0}$, we define
\[
\mu_{q}^{b,d} := \begin{cases}
    \lim_{\epsilon \to 0^+} \left( \left(\beta_q^{b +\epsilon, d - \epsilon} - \beta_q^{b - \epsilon, d - \epsilon}\right) - \left( \beta_q^{b +\epsilon, d + \epsilon} - \beta_q^{b -\epsilon, d + \epsilon} \right) \right)& \textrm{if } d < \infty, \\
    \lim_{\epsilon \to 0^+} \left(\beta_q^{b +\epsilon, \infty} - \beta_q^{b - \epsilon, \infty}\right)& \textrm{otherwise},
\end{cases} 
\] if the index set $T$ is not discrete. For discrete $T$, we define \[\mu_{q}^{b,d} := \begin{cases}
\left(\beta_q^{b',d} - \beta_q^{b,d}\right) - \left(\beta_q^{b',d'} - \beta_q^{b,d'}\right)& \textrm{if } d < \infty,\\
\beta_q^{b',\infty} - \beta_q^{b,\infty}& \textrm{otherwise,}
\end{cases}\]
where $d' = \sup\{\tau \in T \mid \tau < d\}$ and $b' = \sup\{\tau \in T \mid \tau < b\}$. 

The first difference on the right-hand side counts the independent $p$-cycles born at $K_b$ that are alive until $K_d$, while the second difference counts those born at $K_b$ that are still alive after $K_d$. Plotting each point $(b,d)$ with multiplicity $\mu_q^{b,d}$ yields the $q$-th persistence diagram of the filtration. This coordinate system is referred to as the \emph{birth--death coordinate system}. Each point represents a $q^{\textrm{th}}$ homology cycle, with its vertical distance from the diagonal corresponding to its persistence. Since the multiplicities are defined only for $b < d$, all points lie above the diagonal of $\mathcal{R}$. 

We now introduce a method for converting a persistence diagram (PD) into a vector representation that preserves interpretability with respect to the original diagram. This vector representation is known as a \emph{persistence image} (PI)~\cite{PI17}.

Let \( B \) be a persistence diagram. Define a linear transformation \( M : \mathcal{R} \to \mathcal{R} \) by \( M(x, y) = (x, y - x) \) when \( y \neq \infty \), and \( M(x, y) = (x, y) \) when \( y = \infty \). We denote by \( M(B) \) the resulting multiset in \emph{birth--persistence coordinates}, obtained by mapping each point \( (x, y) \in B \) to its image \( M(x,y) \in M(B) \).

Let $\phi_u:\mathcal{R} \to \R$ be an extended Gaussian kernel centered at $u = (u_x, u_y) \in \mathcal{R}$ and variance $\sigma^2$ defined by

\begin{equation} \label{eq: phi_u}
    \phi_u(x, y) = \begin{cases}
        \frac{1}{2\pi\sigma^2} e^{-[(x - u_x)^2 + (y - u_y)^2] / 2\sigma^2}, & \text{ if } u_y\neq \infty, \\
        \frac{1}{2\pi\sigma^2}\delta_{y, \infty}e^{-[(x - u_x)^2] / 2\sigma^2}, & \text{ otherwise},
    \end{cases}
\end{equation}
where $\delta$ is the Dirac-Delta function.

Next, fix a nonnegative uniformly bounded weighting function \( f : \mathcal R \to \R \) that is zero on the horizontal axis. On $\R^2$, $f$ is continuous and piecewise differentiable. The requirement that \( f \) vanishes on the horizontal axis ensures that points near the diagonal, which correspond to topological noise, are disregarded. Using these, we can transform a PD into a scalar function over the extended plane.

\begin{definition}\label{def:ps}
    For a persistence diagram \( B \), the associated \emph{persistence surface} \( \rho_B : \mathcal{R} \to \R \) is defined as
    \[
    \rho_B(z) = \sum_{u \in M(B)} f(u) \, \phi_u(z).
    \]
\end{definition}

To obtain a finite-dimensional vector representation, we discretize a bounded region of the plane and integrate \( \rho_B \) over each cell of the discretization. Specifically, we fix a grid consisting of \( n \) cells (pixels) and assign to each pixel the integral of \( \rho_B \) over that region.

\begin{definition}\label{def:pi}
    Let \( Z \subseteq \mathcal{R} \) be a bounded grid with \( n \) pixels, and let \( P \) denote one such pixel. For a persistence diagram \( B \), the corresponding \emph{persistence image} is the collection of values
    \[
    I_{\rho_B}(P) = \iint_P \rho_B(y, z) \, dy \, dz.
    \]
\end{definition}

Persistence images provide an effective means of combining PDs across different homological dimensions into a unified representation. For instance, let persistence diagrams for \( H_0, H_1, \ldots, H_k \) are computed in an experiment. One can concatenate the corresponding PI vectors for \( H_0, H_1, \ldots, H_k \) into a single feature vector that encodes all homological dimensions simultaneously, which can then be used as input for machine learning algorithms.

\subsection{Filtrations}\label{sec:filtration}

To apply the persistent Laplacian framework to discrete data, such as point clouds or graphs, one must first construct a simplicial filtration $K = \{ K_t \}_{t \in T}$ that captures the geometry or topology of the data at varying scales.

\begin{definition}\label{def:filtration}

A \emph{filtration} on a simplicial complex is a family of nested subcomplexes
\[
K = \{ K_t \}_{t \in T}, \qquad \text{with } s \le t \implies K_s \subseteq K_t,
\]
indexed by a totally ordered set \(T \subseteq \mathbb{R}\).  
Equivalently, a filtration may be generated by a function \(f : K \to \mathbb{R}\) by
\[
K_t = f^{-1}((-\infty, t]).
\]    
\end{definition}

Filtrations encode multiscale structure and serve as the underlying domain for both persistent homology and persistent Laplacians.

\begin{figure}[t]
    \centering
    \includegraphics[width=0.9\linewidth]{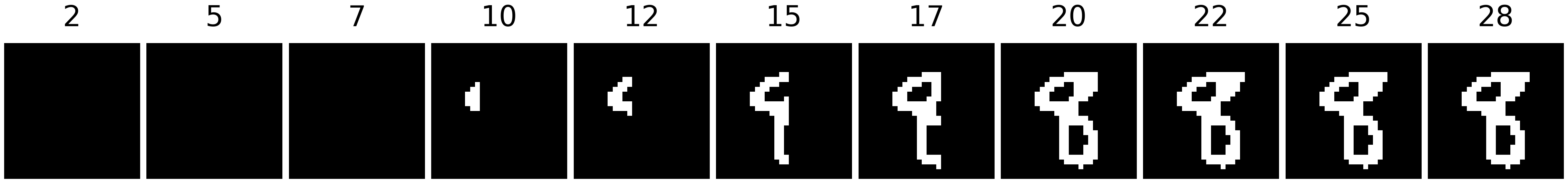}
    \caption{Height$(\to)$ filtration of a figure 8 MNIST sample, shown from left to right as the filtration value increases.}
    \label{fig:mnist-height-filtration}
\end{figure}

\section{Signatures}
\label{sec: signatures}

Our goal is to develop a variant of a persistence diagram (PD) that incorporates information from the persistent Laplacian (PL) through an additional variable, which we call a \emph{signature}. We refer to this variant as a persistent Laplacian diagram (PLD). Examples of signatures and their computation are presented in~\Cref{subsec:sign}.

\begin{definition}[Signatures] \label{def:signatures}
Let $\mathcal{S}_n$ denote the space of real symmetric matrices of finite size $n$.  
A \emph{signature} is any map
\[
s : \mathcal{S}_n \longrightarrow \R
\]
that assigns a real number to each such matrix of size $n$ satisfying $s(Q^{-1} A Q)=s(A)$ for every invertible $Q$. For a subspace $\mathcal A \subset \mathcal S_n$, we call the signature $s$ is \emph{admissible} over $\mathcal A$ when $\sup_{\mathcal A}\|s\|_\infty < \mathcal{M}$ for some constant $\mathcal{M}$ depending on $n$. Unless stated otherwise, we take $\mathcal A$ to be $\mathcal S_n$.
\end{definition}

\begin{remark}
Admissibility is a condition relative to the chosen data class and signature, rather than a universal property of all persistent Laplacians of arbitrary size. It may fail for adversarial families, such as complexes of unbounded size. In practical finite-data settings, however, the relevant family of matrices is bounded, so standard signatures such as the trace, spectral entropy, and geometric profiles remain uniformly bounded on the observed data class. Extending admissibility to settings involving complexes of unbounded size is an important direction for future work.
\end{remark}

\subsection{Persistent Laplacian Diagram}
Intuitively, persistent homology tracks the birth and death of homology classes along a filtration. Each off--diagonal point \((b,d)\) in the persistence diagram records an interval during which a nontrivial \(q\)-cycle exists before merging into another generator. In contrast, the persistent Laplacian can be viewed as a tool that, at the moments where generators are created or annihilated, describes how these events occur in the underlying chain complex. 

Concretely, suppose new \(q\)-generators are born at filtration values \(2,3,5\) and never die, so persistence diagram records \((2,\infty),(3,\infty),(5,\infty)\). Persistent Laplacians, by examining operators such as \(\Delta_q^{3,2}\), \(\Delta_q^{5,2}\), and \(\Delta_q^{5,3}\), encode not only when generators appear, but also additional information about their evolution under filtration.

Fix a homological degree \(q\ge 0\) and a filtration \(\mathbf{K} = \{K_t\}_{t\in T}\) as in Section~\ref{sec:filtration} with $K_\infty := \bigcup_{t\in T}K_t$. 
Let \(\mu_q^{b,d}\) denote the multiplicity of the point \((b,d)\) in the \(q\)-th persistence diagram $B_q$, and  \[
\mathcal{B}_q := \{\, b \in T : \exists d\in T\cup\{\infty\} \text{ with } \mu_q^{b,d} \neq 0 \,\},
\quad
\mathcal{D}_q := \{\, d \in T : \exists b\in T \text{ with } \mu_q^{b,d} \neq 0 \,\}.
\]
We define \(C_{B_q}\) the set of tuples as
\[
C_{B_q}:= \left\{\, (b,d) \in T \times \left(T \cup \{\infty\} \right) : d>b ~\textrm{ and }~ b,d \in \mathcal{B}_q \cup \mathcal{D}_q\right\}.
\]



For each pair $(b,d)\in C_{B_q}$, the filtration provides a simplicial pair $
K_b \hookrightarrow K_d$,
and hence a persistent Laplacian
$\Delta^{K_d,K_b}_q : C_q^{K_b} \longrightarrow C_q^{K_b},
$ as in Section~\ref{sec:prelim}.  
The operator $\Delta^{K_d,K_b}_q$ is represented, with respect to a fixed basis of \(C_q^{K_b}\), by a real symmetric matrix of size \(n^{K_b}_q \), where \(n^{K_b}_q = \dim C_q^{K_b}\).

Recall that $\mathcal{R} = \R \times \overline{\R}$ is the domain of $(b, d)$ tuples in the $q$-th persistence diagram $B_q$. The \emph{\(q\)-th persistent Laplacian diagram} (PLD) is \(\Gamma(S_{B_q}) = \left\{(u,S_{B_q}(u)) : u \in \mathcal R\right\}\) the graph of function $S_{B_q} : \mathcal R \longrightarrow \R$ defined as
\[
S_{B_q}(u) = \begin{cases}
    s\bigl(\Delta^{K_d,K_b}_q\bigr),& \textrm{if } u \in C_{B_q}, \\
    0,& \textrm{otherwise,}
\end{cases}
\]
where $s$ be a signature. By the abuse of notation, we denote PLD as $C_{B_q}$ and $\Gamma(S_{B_q})$ interchangeably if $s$ is fixed. 

\begin{remark}
A PLD is a signature-dependent and the underlying operators cannot in general be recovered from one scalar value per pair.
\end{remark}

We extend the notion of the \( p \)-Wasserstein distance, originally defined between two PDs, to define an analogous metric between two PLDs. Recall that the infinity norm is given by\(\| f - g \|_{\infty} = \sup_{x \in X} |f(x) - g(x)|\) for functions \( f, g : X \to \R \) and \(\|x\|=\sup_{1\leq i\leq n}|x_i|\) for a vector $x=(x_1,\ldots,x_n)\in\R^n$. Using this, the \( p \)-Wasserstein distance between two PLDs \( C_{B_q} \) and \( C_{B'_q} \) is defined by
\begin{equation} \label{eq:wasserstein}
W_p(C_{B_q}, C_{B'_q}) = \inf_{\gamma : C_{B_q} \rightarrow C_{B'_q}} \left( \sum_{u \in C_{B_q}} \| u - \gamma(u) \|_{\infty}^p \right)^{1/p},
\end{equation}
where \( 1 \leq p < \infty \) and \( \gamma \) ranges over all standard maximum matching between \( C_{B_q} \) and \( C_{B'_q} \), where the surplus points of diagram corresponds to the closest point on the diagonal~\cite{CEH07}.

For the case \( p = \infty \), we obtain the standard \emph{bottleneck distance}, defined analogously by
\[
W_{\infty}(C_{B_q}, C_{B'_q}) = \inf_{\gamma : C_{B_q} \rightarrow C_{B'_q}} \left( \sup_{u \in C_{B_q}} \| u - \gamma(u) \|_{\infty}\right).
\]

These metrics provide a quantitative measure of the similarity between the homological features of two data sets, as captured by their corresponding persistence diagrams or persistent Laplacian diagrams.

\subsection{Persistent Laplacian Image}

Given a \(q\)-th persistence diagram $B_q$, we now define its associated \emph{persistent Laplacian surface and image} in direct analogy with the persistence surface and persistence image of~\cite{PI17}, but using the values of the signature as amplitudes.

We use the extended Gaussian kernel $\phi_u:\mathcal{R} \to \R$ from \Cref{eq: phi_u} and the weighting function $f:\mathcal{R}\to\R_{\ge 0}$. Recall that $f$ is nonnegative, vanishes on the diagonal, and is continuous and piecewise differentiable on $\R^2$. Then we can define persistent Laplacian surface and persistent Laplacian image similarly as~\Cref{def:ps} and~\Cref{def:pi}. Recall that $M : \mathcal R \to \mathcal R$ by $M(x,y) = (x,y-x)$.

Let \(C_{B_q}\) be the support of function \(S_{B_q}\) as above. 
\begin{definition}[Persistent Laplacian surface and image]
The \emph{persistent Laplacian surface} associated with persistent Laplacian diagram \(C_{B_q}\) is the function
\[
\rho_{C_{B_q}} : \mathcal R\longrightarrow \R,
\qquad
\rho_{C_{B_q}}(x,y)
    = \sum_{u \in C_{B_q}} f\circ M(u)\, S_{B_q}(u)\, \phi_{M(u)}(x,y).
\]

Let \(Z \subset \mathcal R\) be a bounded grid decomposed into pixels \(P\).  The corresponding persistent Laplacian image (PLI) is the collection of values
\[
I_{C_{B_q}}(P)
    = \iint_P \rho_{C_{B_q}}(x,y)\,dx\,dy,
\]
one scalar per pixel \(P\).  Flattening these pixel values yields a finite--dimensional feature vector associated with the \(q\)-th persistent Laplacian.
\end{definition}

\begin{figure}[t]
    \centering
    \begin{subfigure}[t]{0.48\linewidth}
        \centering
        \includegraphics[width=\linewidth]{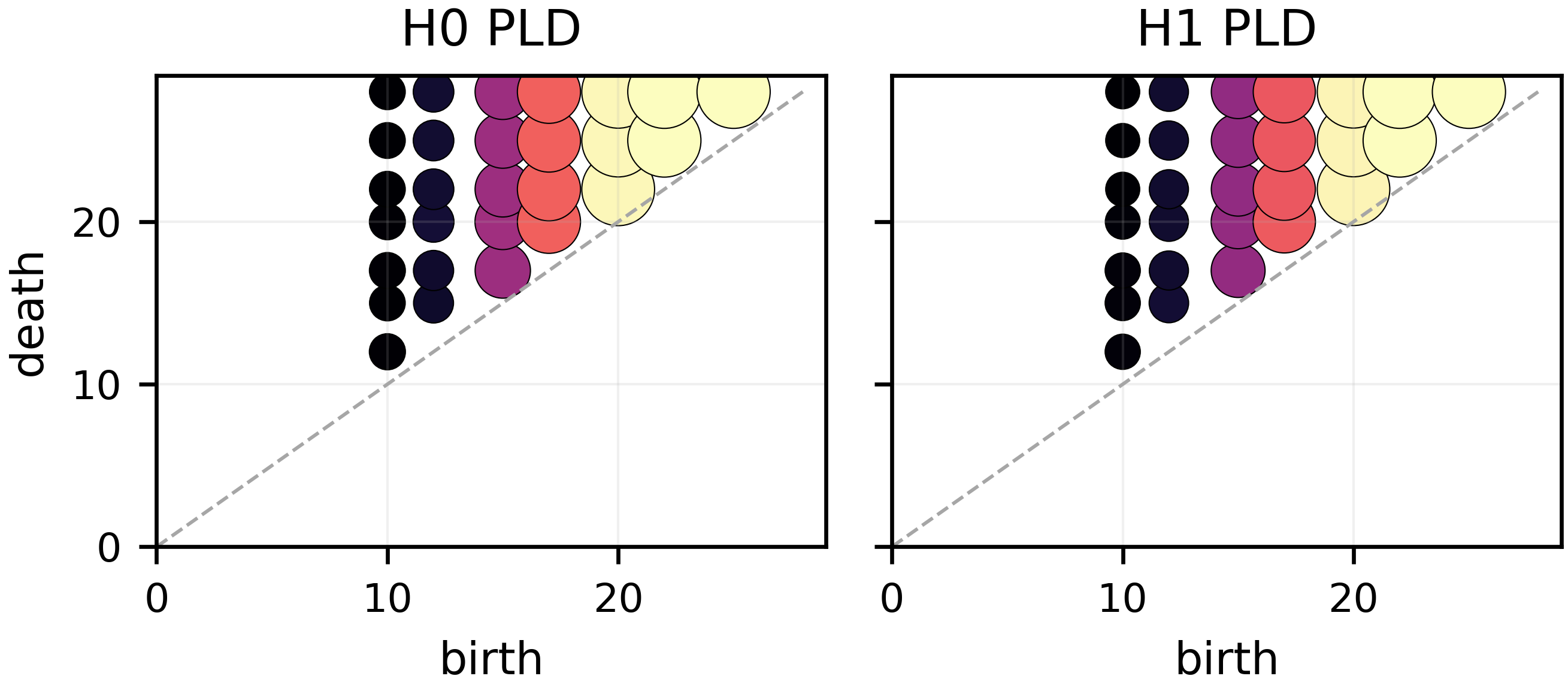}
        \caption{Persistent Laplacian diagram (PLD).}
        \label{fig:sample-pld}
    \end{subfigure}
    \hfill
    \begin{subfigure}[t]{0.48\linewidth}
        \centering
        \includegraphics[width=\linewidth]{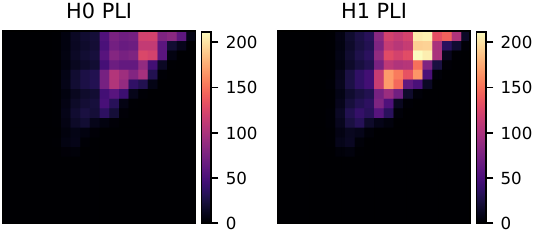}
        \caption{Persistent Laplacian image (PLI).}
        \label{fig:sample-pli}
    \end{subfigure}
    \caption{Persistent Laplacian diagram (left) and persistent Laplacian image (right) constructed from the MNIST height filtration shown in Figure~\ref{fig:mnist-height-filtration}. The circles in PLD represent the relative values of the chosen persistent Laplacian signature(trace) across birth--death pairs induced by the filtration. The PLI is obtained by smoothing and discretizing this diagram into a vectorized representation.}
    \label{fig:sample-pld-pli}
\end{figure}

\section{Stability}
\label{sec: stability}

We study stability directly at the level of persistent Laplacian diagrams. This was discussed for the stability of persistent diagrams in~\cite{PI17}. The admissibility condition in Definition~\ref{def:signatures} provides a uniform bound on the signature amplitudes, while the geometric displacement of PLD points is measured by their Wasserstein distance. For Gaussian kernels, the persistent Laplacian image satisfies the following stability bound. Let $\mathcal M_S:=\max\{\|S_{B_q}\|_\infty,\|S_{B'_q}\|_\infty\}$.

\begin{restatable}{theorem}{GaussianImageStability}\label{thm:Gauss}
The persistent Laplacian image associated with a Gaussian kernel is stable with respect to the $1$-Wasserstein distance between the PLDs $C_{B_q}$ and $C_{B'_q}$. More precisely, for $p=1,2,$ and $\infty$,
\[
\|I_{\rho_{C_{B_q}}}-I_{\rho_{C_{B'_q}}}\|_p
\leq \mathcal M_S
\left(\sqrt{5}|\nabla f|+\sqrt{\frac{10}{\pi}}\frac{\|f\|_\infty}{\sigma}\right)
W_1(C_{B_q},C_{B'_q}).
\]
\end{restatable}

The proof is provided in Appendix~\ref{sec:proofs}.

\section{Experiments}
\label{sec:experiments}

The \textbf{MNIST} dataset is a standard benchmark consisting of $28 \times 28$ grayscale images of handwritten digits \cite{lecun2010mnist}. In our experiments, we follow the same height filtration used in the baseline of \cite{davies2023persistent}, proceeding from left to right across the image.

Another experiment is on the \textbf{QM7} \cite{blum, rupp}, which contains molecular data in the Coulomb matrix $M=[m_{i,j}]_{i,j\in I}$ for a given set of nuclei $I$. In this matrix, the entries $m_{i,j}$ represent the electrostatic interactions between atomic nuclei within the molecule \cite{rupp}. To extract topological features from the \textbf{QM7} dataset, we adopt a distance-dependent complex generation framework. We compute the pairwise atomic distances derived from the Coulomb matrix elements, $C_{ij} = q_i q_j / \|\mathbf{r}_i - \mathbf{r}_j\|$, where $q_i, q_j$ represent atomic charges and $\mathbf{r}_i, \mathbf{r}_j$ are their coordinates. Our construction follows the filtration 2 in the baseline experiment protocol of \cite{davies2023persistent}. If two vertices apart from each other with distance $d$, the constituent 0-simplices (atoms) are admitted into the 0-degree complex when the filtration value exceeds $d$. Consequently, the vertex set grows dynamically throughout the filtration rather than being entirely present from the initiation.

The reproduced PL baseline of~\cite{davies2023persistent} extracts PL-based features by padding and truncating the eigenvalue sequence of the persistent Laplacian, whereas our method instead applies the proposed persistent Laplacian image (PLI) construction to obtain a fixed-dimensional representation. In both cases, features are extracted from the $0$- and $1$-dimensional persistent Laplacians and then concatenated into a single flat feature vector. We additionally evaluate the standard persistence image (PI) of~\cite{PI17}, formed separately from the $H_0$ and $H_1$ persistence diagrams and concatenated into a single feature vector.

For Table~\ref{tab:mnist_accuracy_metrics} and Table~\ref{tab:persistent_laplacian_metrics}, we use filtration resolution \(r_f=12\) for MNIST and \(r_f\in\{8,12\}\) for QM7. For PI and PLI, we evaluate image resolutions \(r_i\in\{10,20,30\}\). We compare the reproduced raw persistent Laplacian (PL) baseline, the ordinary persistence image (PI), and PLI using the trace and spectral-entropy signatures.

For MNIST, we report the mean test accuracy and its sample standard deviation across these five runs. For QM7, we report the mean and sample standard deviation across the five fold-level averages, corresponding to 15 model fits per configuration.

\section{Results and Discussions}
\label{sec: results}

\begin{table}[t]
\centering
\caption{
MNIST classification results at filtration resolution \(r_f=12\) in accuracy.
}
\label{tab:mnist_accuracy_metrics}
\small
\setlength{\tabcolsep}{5pt}
\begin{tabular}{@{}l|c|c|c@{}}
\toprule
\textbf{Representation} & \(\boldsymbol{r_i=10}\) &
\(\boldsymbol{r_i=20}\) & \(\boldsymbol{r_i=30}\) \\
\midrule
PI & \(0.6566 \pm 0.0009\) & \(0.6564 \pm 0.0009\) & \(0.6561 \pm 0.0022\) \\
PLI (trace) & \(0.8772 \pm 0.0015\) &
\(\mathbf{0.8834 \pm 0.0010}\) & \(0.8821 \pm 0.0011\) \\
PLI (spectral entropy) & \(0.8545 \pm 0.0017\) &
\(0.8611 \pm 0.0026\) & \(0.8601 \pm 0.0031\) \\
\midrule
PL baseline & \multicolumn{3}{c}{\(0.8328 \pm 0.0018\)} \\
\bottomrule
\end{tabular}
\end{table}

\begin{table}[t]
\centering
\caption{
QM7 regression results in MAE and RMSE.}
\label{tab:persistent_laplacian_metrics}
\small
\setlength{\tabcolsep}{3.5pt}
\begin{tabular}{@{}lcccc@{}}
\toprule
\textbf{Representation} &
\(\boldsymbol{r_f}\) &\textbf{Image Resolution}
\(\boldsymbol{r_i}\) &
\textbf{MAE}  &
\textbf{RMSE}\\
\midrule
PL baseline
 & 12 & -- & \(50.676 \pm 0.673\) & \(64.849 \pm 1.112\) \\
\midrule
PI
 & 12 & 10 & \(130.514 \pm 1.281\) & \(165.693 \pm 1.689\) \\
 & 12 & 20 & \(130.898 \pm 1.325\) & \(166.074 \pm 1.694\) \\
 & 12 & 30 & \(130.542 \pm 0.941\) & \(165.730 \pm 1.146\) \\
\midrule
PLI (trace)
 &  8 & 10 & \(37.789 \pm 1.065\) & \(48.528 \pm 1.713\) \\
 &  8 & 20 & \(37.563 \pm 0.943\) & \(47.754 \pm 1.390\) \\
 &  8 & 30 & \(37.807 \pm 0.908\) & \(47.985 \pm 1.407\) \\
 & 12 & 10 & \(37.668 \pm 0.764\) & \(49.847 \pm 1.813\) \\
 & 12 & 20 & \(36.844 \pm 1.016\) & \(47.987 \pm 1.844\) \\
 & 12 & 30 &
   \(\mathbf{36.722 \pm 0.813}\) &
   \(\mathbf{47.732 \pm 1.424}\) \\
\midrule
PLI (spectral entropy)
 &  8 & 10 & \(39.192 \pm 0.761\) & \(50.131 \pm 1.467\) \\
 &  8 & 20 & \(39.268 \pm 0.683\) & \(50.439 \pm 1.554\) \\
 &  8 & 30 & \(39.247 \pm 0.687\) & \(50.345 \pm 1.284\) \\
 & 12 & 10 & \(39.297 \pm 1.117\) & \(52.712 \pm 3.437\) \\
 & 12 & 20 & \(40.007 \pm 0.554\) & \(53.422 \pm 2.507\) \\
 & 12 & 30 & \(40.410 \pm 0.544\) & \(53.941 \pm 2.880\) \\
\bottomrule
\end{tabular}
\end{table}

\begin{table}[t]
\centering
\small
\caption{Published comparison methods from Ref.~\cite{davies2023persistent}.}
\label{tab:davies_published}
\begin{tabular}{lcccc}
\toprule
\textbf{Task / metric} & \textbf{PH} & \textbf{GL} & \textbf{CL} & \textbf{PL} \\
\midrule
MNIST / accuracy & $0.625\pm0.010$ & $0.768\pm0.006$ & $0.798\pm0.006$ & $0.821\pm0.011$ \\
QM7-CM Filtration 2 / MAE & $45.47\pm1.438$ & $23.19\pm0.809$ & $23.33\pm0.785$ & $21.94\pm0.737$ \\
\bottomrule
\end{tabular}
\end{table}

Tables~\ref{tab:mnist_accuracy_metrics} and
\ref{tab:persistent_laplacian_metrics} summarize the controlled
evaluation. On MNIST, trace-based PLI achieves its highest mean
accuracy at \(r_i=20\), obtaining approximately \(0.8834\), compared with \(0.8328\) for the reproduced PL baseline. This
corresponds to an improvement of approximately \(5.1\) percentage
points. Spectral-entropy PLI also exceeds the PL baseline for all
three tested image resolutions, whereas the ordinary PI performs
substantially worse under the tested construction.

On QM7, trace-based PLI at \((r_f,r_i)=(12,30)\) yields the lowest
mean errors, with an MAE of \(36.722\pm0.813\) and an RMSE of
\(47.732\pm1.424\). Relative to the reproduced PL baseline, these
values represent reductions of approximately \(27.5\%\) in MAE and
\(26.4\%\) in RMSE. The ordinary PI is not competitive on QM7 under the present filtration.

These results suggest that PLI provides a flexible fixed-dimensional representation of PL, where the choice of signature plays the role of a tunable design parameter. We also provide the experiment results by varing the signatures to find representations that achieve better predictive performance, see Table~\ref{tab:old_mnist} and Table~\ref{tab:old_qm7}. For broader TDA context, Table~\ref{tab:davies_published} also lists PH, graph-Laplacian (GL), combinatorial-Laplacian (CL), and persistent-Laplacian (PL) results reported by Ref.~\cite{davies2023persistent}.

\section{Conclusions}
\label{sec:conclusions}

In this work, we introduced a signature-based vectorization framework for persistent Laplacians through persistent Laplacian diagrams and persistent Laplacian images. Our construction turns persistent Laplacians into fixed-dimensional representations that can be used directly in downstream learning tasks, while remaining flexible in the choice of spectral signature. On the theoretical side, we established stability guarantees for Gaussian persistent Laplacian images with respect to Wasserstein perturbations at the PLD level. On the empirical side, experiments on MNIST and QM7 showed that this representation is effective in practice and that the choice of signature has a substantial impact on performance, with the trace emerging as a particularly strong option among those considered.

Potential societal impact should be understood primarily through the application domains in which topological data analysis and persistent Laplacians have already shown practical value. Prior work on persistent Laplacians, particularly in molecular and biomolecular settings, suggests that they can capture both topological and spectral information relevant to scientific prediction tasks. At present, we view the main impact of this work as methodological: providing a stable and flexible representation of persistent Laplacians for downstream learning tasks.

Despite these strengths, the present work has several limitations. First, while the PLI framework is flexible, its practical performance depends on design choices such as the signature, the filtration resolution, and the image resolution. Our experiments show that these choices matter substantially---for example, the trace performs best among the signatures we considered---but we do not yet provide a principled, data-adaptive rule for selecting them across tasks. Second, the computational cost can be significant, since constructing PLIs requires evaluating persistent Laplacians over many birth--death pairs and repeatedly computing spectral information, especially when fine filtration sampling is needed. One way to mitigate this cost is to develop signatures based on only a small subset of the spectrum, since computing only a few small eigenvalues by iterative methods can reduce computational cost for spectrum substantially while preserving much of the predictive value of the representation. Appendix~\ref{app:pl-algorithms} gives the explicit construction together with time and memory bounds. Third, our empirical evaluation is limited to MNIST and QM7. Additional experiments on other datasets would further strengthen the empirical evidence for the generality of the proposed framework.

Several additional directions remain for future research. One is  to evaluate the framework on a wider range of datasets and application domains, including benchmarks such as PDBind, to better understand where persistent Laplacian images are most effective. On the theoretical side, a more detailed analysis of when similar birth--death structure gives rise to similar persistent Laplacian signatures may help clarify the information added by PL beyond persistence diagrams alone. Another direction is to identify sharper and more practically useful conditions under which signatures satisfy the stability requirements of the framework.

\newpage
\bibliographystyle{plain}
\bibliography{bib}


\appendix

\section{Proofs}\label{sec:proofs}

Let $C_{B_q}$ and $C_{B'_q}$ be two PLDs. Let $\gamma$ denote a Wasserstein matching in the sense of~\eqref{eq:wasserstein}; thus, surplus points are paired with the closest diagonal. Let $\mathcal M_S:=\max\{\|S_{B_q}\|_\infty,\|S_{B'_q}\|_\infty\}$.
Throughout this section, we work on the grid $Z\subset\R^2$ and restrict attention to $\R^2$ rather than the larger domain $\mathcal R$.

For a differentiable function \( h : \R^2 \to \R \), let \( |\nabla h| \) denote the maximal norm of the gradient vector of \( h \), that is, the largest directional derivative of \( h \). By the fundamental theorem of calculus for line integrals, for all \( u, v \in \R^2 \), we have
\[
|h(u) - h(v)| \leq |\nabla h| \, \|u - v\|_2.
\]

For a differentiable probability distribution \( \phi_u \) with mean \( u = (u_x, u_y) \in \R^2 \), we may safely denote \( |\nabla \phi_u| \) by \( |\nabla \phi| \) and \( \|\phi_u\|_{\infty} \) by \( \|\phi\|_{\infty} \), since both the maximal directional derivative and the supremum of a differentiable probability distribution are invariant under translation.  
For any \( z \in \R^2 \), we then have
\[
|\phi_u(z) - \phi_v(z)| = |\phi_u(z) - \phi_u(z + u - v)| \leq |\nabla \phi| \, \|u - v\|_2,
\]
which implies
\[
\|\phi_u - \phi_v\|_{\infty} \leq |\nabla \phi| \, \|u - v\|_2.
\]

\begin{lemma}\label{lem:wf}
Let \( f \) be a weighting function. For all \( u, v \in \R^2 \), we have
\[
\| f(u)\phi_u - f(v)\phi_v \|_{\infty} \leq \left( \|f\|_{\infty} |\nabla \phi| + \|\phi\|_{\infty} |\nabla f| \right) \|u - v\|_2.
\]
\end{lemma}

\begin{proof}
Using the triangle inequality and the pointwise bound for the directional derivative, we obtain for any $z\in\mathbb{R}^2$
\begin{align*}
    |f(u)\phi_u(z)-f(v)\phi_v(z)|
    &\le |f(u)|\;|\phi_u(z)-\phi_v(z)|+|\phi_v(z)|\;|f(u)-f(v)| \\
    &\le \|f\|_{\infty}\,|\phi_u(z)-\phi_v(z)|+\|\phi\|_{\infty}\,|f(u)-f(v)|.
\end{align*}
Since $\phi$ and $f$ are differentiable, the mean value inequality (or the fundamental theorem of calculus for line integrals) gives
\[
|\phi_u(z)-\phi_v(z)| \le |\nabla\phi|\,\|u-v\|_2,\qquad
|f(u)-f(v)| \le |\nabla f|\,\|u-v\|_2,
\]
where $|\nabla\phi|$ and $|\nabla f|$ denote the supremum norms of the gradients of $\phi$ and $f$, respectively. Combining these estimates and taking the supremum over $z$ yields
\[
\|f(u)\phi_u-f(v)\phi_v\|_{\infty}
\le \big(\|f\|_{\infty}|\nabla\phi|+\|\phi\|_{\infty}|\nabla f|\big)\|u-v\|_2,
\]
as claimed.
\end{proof}

\begin{theorem}\label{thm:pls}
The persistent Laplacian surface \( \rho \) is stable with respect to the \( 1 \)-Wasserstein distance between the PLDs $C_{B_q}$ and $C_{B'_q}$:
\[
\|\rho_{C_{B_q}}-\rho_{C_{B'_q}}\|_{\infty}
\le 2\sqrt{2}\mathcal M_S\big(\|f\|_{\infty}|\nabla\phi|+\|\phi\|_{\infty}|\nabla f|\big)\,W_1(C_{B_q},C_{B'_q}).
\]
\end{theorem}

\begin{proof}
Let $\gamma$ be an optimal $1$-Wasserstein matching between $C_{B_q}$ and $C_{B'_q}$. Using the zero extension on diagonal matches, we have
    \begin{align*}
    \|\rho_{C_{B_q}}-\rho_{C_{B'_q}}\|_{\infty}
    &\le \sum_{(u,v)\in\gamma}
    \left\|f\circ M(u)S_{B_q}(u)\phi_{M(u)}
    -f\circ M(v)S_{B'_q}(v)\phi_{M(v)}\right\|_\infty\\
    &\le \sqrt{2}\mathcal M_S
    \big(\|f\|_{\infty}|\nabla\phi|+\|\phi\|_{\infty}|\nabla f|\big)
    \sum_{(u,v)\in\gamma}\|M(u)-M(v)\|_\infty.
    \end{align*}
The pairwise estimate follows from Lemma~\ref{lem:wf}, and
$\|M(u)-M(v)\|_{\infty}\le 2\|u-v\|_{\infty}$. Therefore,
\[
\|\rho_{C_{B_q}}-\rho_{C_{B'_q}}\|_{\infty}
\le 2\sqrt{2}\mathcal M_S\big(\|f\|_{\infty}|\nabla\phi|+\|\phi\|_{\infty}|\nabla f|\big)\,W_1(C_{B_q},C_{B'_q}).
\]
\end{proof}

\begin{theorem}\label{thm:pls2}
The persistence image $I_{\rho_{C_{B_q}}}$ is stable with respect to the $1$-Wasserstein distance between the PLDs. More precisely, if $A$ is the maximum area of any pixel, $A'$ is the total image area, and $n$ is the number of pixels, then
\[
\|I_{\rho_{C_{B_q}}}-I_{\rho_{C_{B'_q}}}\|_\infty
\leq 2\sqrt{2}A\mathcal M_S\big(\|f\|_{\infty}|\nabla\phi|+\|\phi\|_{\infty}|\nabla f|\big)W_1(C_{B_q},C_{B'_q}),
\]
\[
\|I_{\rho_{C_{B_q}}}-I_{\rho_{C_{B'_q}}}\|_1
\leq 2\sqrt{2}A'\mathcal M_S\big(\|f\|_{\infty}|\nabla\phi|+\|\phi\|_{\infty}|\nabla f|\big)W_1(C_{B_q},C_{B'_q}),
\]
and
\[
\|I_{\rho_{C_{B_q}}}-I_{\rho_{C_{B'_q}}}\|_2
\leq 2\sqrt{2n}A\mathcal M_S\big(\|f\|_{\infty}|\nabla\phi|+\|\phi\|_{\infty}|\nabla f|\big)W_1(C_{B_q},C_{B'_q}).
\]
\end{theorem}

\begin{proof}
For any pixel $P$ with area $A(P)$,
    \begin{align*}
\left|I_{\rho_{C_{B_q}}}(P)-I_{\rho_{C_{B'_q}}}(P)\right|
&\leq A(P)\|\rho_{C_{B_q}}-\rho_{C_{B'_q}}\|_\infty.
    \end{align*}
Applying Theorem~\ref{thm:pls}, summing over pixels, and using
$\|x\|_2\leq\sqrt n\|x\|_\infty$ gives the three bounds.
\end{proof}

\begin{theorem}\label{thm:GaussSurface}
The persistent Laplacian surface associated with a Gaussian kernel is stable with respect to the $1$-Wasserstein distance between $C_{B_q}$ and $C_{B'_q}$:
\[
\|\rho_{C_{B_q}}-\rho_{C_{B'_q}}\|_1
\leq \mathcal M_S
\left(\sqrt{5}|\nabla f|+\sqrt{\tfrac{10}{\pi}}\tfrac{\|f\|_\infty}{\sigma}\right)
W_1(C_{B_q},C_{B'_q}).
\]
\end{theorem}

\begin{proof}
Let $\gamma$ be an optimal $1$-Wasserstein matching between the two PLDs. Using the zero extension for diagonal matches,
    \begin{align*}
\|\rho_{C_{B_q}}-\rho_{C_{B'_q}}\|_1
&\leq \sum_{(u,v)\in\gamma}
\left\|f\circ M(u)S_{B_q}(u)g_{M(u)}
-f\circ M(v)S_{B'_q}(v)g_{M(v)}\right\|_1\\
&\leq \mathcal M_S
\left(|\nabla f|+\sqrt{\tfrac{2}{\pi}}\tfrac{\|f\|_\infty}{\sigma}\right)
\sum_{(u,v)\in\gamma}\|M(u)-M(v)\|_2\\
&\leq \mathcal M_S
\left(\sqrt{5}|\nabla f|+\sqrt{\tfrac{10}{\pi}}\tfrac{\|f\|_\infty}{\sigma}\right)
\sum_{(u,v)\in\gamma}\|u-v\|_\infty\\
&=\mathcal M_S
\left(\sqrt{5}|\nabla f|+\sqrt{\tfrac{10}{\pi}}\tfrac{\|f\|_\infty}{\sigma}\right)
W_1(C_{B_q},C_{B'_q}).
    \end{align*}
The second inequality follows from Lemma 3 of~\cite{PI17}, and the third uses
$\|M(u)-M(v)\|_2\leq\sqrt5\|u-v\|_\infty$.
\end{proof}

\GaussianImageStability*

\begin{proof}
For the $L_1$ norm of the image vector,
    \begin{align*}
\|I_{\rho_{C_{B_q}}}-I_{\rho_{C_{B'_q}}}\|_1
&=\sum_P\left|\iint_P\rho_{C_{B_q}}\,dxdy-\iint_P\rho_{C_{B'_q}}\,dxdy\right|\\
&\leq\iint_{\mathbb R^2}|\rho_{C_{B_q}}-\rho_{C_{B'_q}}|\,dxdy\\
&=\|\rho_{C_{B_q}}-\rho_{C_{B'_q}}\|_1\\
&\leq\mathcal M_S
\left(\sqrt{5}|\nabla f|+\sqrt{\tfrac{10}{\pi}}\tfrac{\|f\|_\infty}{\sigma}\right)
W_1(C_{B_q},C_{B'_q}),
    \end{align*}
by Theorem~\ref{thm:GaussSurface}. The $p=2$ and $p=\infty$ cases follow from
$\|x\|_2\leq\|x\|_1$ and $\|x\|_\infty\leq\|x\|_1$ for finite-dimensional vectors.
\end{proof}

\section{Direct Examples} \label{sec: examples}

One can observe directly that by adding information on graph Laplacian, we can classify more graphs than PH.

\subsection{Distinct Eigenvalues}\label{sec: distinct eigenvalues}

\begin{figure}[h]
    \centering
\begin{tikzpicture}
    \draw (0,0) to (1,.5);
    \draw (0,0) to (0,1);
    \draw (0,1) to (1,.5);
    \draw (1,.5) to (2,.5);
    \draw (2,.5) to (3,1);
    \draw (2,.5) to (3,0);
    \draw (3,0) to (3,1);
    \fill (0,0) circle(.1);
    \fill (1,.5) circle(.1);
    \fill (0,1) circle(.1);
    \fill (2,.5) circle(.1);
    \fill (3,1) circle(.1);
    \fill (3,0) circle(.1);
    \node at (1.5,-1) {$G_1$};

    \draw (5,0) to (5,1);
    \draw (5,0) to (6,0);
    \draw (5,1) to (6,1);
    \draw (6,1) to (6,0);
    \draw (6,1) to (7,1);
    \draw (6,0) to (7,0);
    \draw (7,0) to (7,1);
    \fill (5,0) circle(.1);
    \fill (5,1) circle(.1);
    \fill (6,0) circle(.1);
    \fill (6,1) circle(.1);
    \fill (7,0) circle(.1);
    \fill (7,1) circle(.1);
    \node at (6,-1) {$G_2$};
\end{tikzpicture}
    \caption{Non-isomorphic codegree pair of graphs $G_1$ and $G_2$. They are indistinguishable by PH in degree filtration.}
    \label{fig:eigenval_dist}
\end{figure}
We can see that the graphs $G_1$ and $G_2$ in \Cref{fig:eigenval_dist} have the same degree sequence $(2,2,2,2,3,3)$. Hence they cannot be distinguished by PH on degree filtration, thereby their PIs are the same. Their combinatorial graph Laplacian spectra $\left(\frac{5+\sqrt{17}}{2},3,3,3,\frac{5-\sqrt{17}}{2},0\right)$ and $(5,3,3,2,1,0)$ are different, so we can distinguish them using eigenvalue based signatures, such as the second smallest eigenvalue $s_{gap}$ of the graph laplacian, which corresponds to the spectrum of Persistant Laplacians at infinity $\Delta^{\infty,3}_0$.

\begin{figure}[h]
    \centering
\begin{tikzpicture}[scale=1]

\begin{scope}[shift={(0,0)}]
    \fill (0,0) circle (.1);
    \node[below](G3) at (0,-0.3) {$G_3 = H_3$ };
\end{scope}

\begin{scope}[shift={(4.5,2)}]
    \draw (0,0) to (-1,.5);      
    \draw (0,0) to (-1,-0.5);    
    \draw (0,0) to (1,.5);       
    \draw (0,0) to (1,-0.5);    

    \draw (-1, .5) to (1, .5);  
    \draw (-1, .5) to (-1, -0.5);

    \draw (1, .5) to (2, 0);     
    \draw (2, 0) to (1, -0.5);   
    \draw (1, -0.5) to (-1, -0.5);

    \fill (0,0) circle(.1);
    \fill (1,.5) circle(.1);
    \fill (-1,.5) circle(.1);
    \fill (-1,-0.5) circle(.1);
    \fill (2,0) circle(.1);
    \fill (1,-0.5) circle(.1);

    \node[align=center] (G4) at (0,-1) {$G_{4}$\\ 4 cycles};
\end{scope}

\begin{scope}[shift={(4.5,-2)}]
    
    \draw (-.5,0)   to (-1,.5);    
    \draw (-.5,0)   to (-1,-0.5);  
    \draw (-.5,0)   to (.5,0);       

    \draw (.5, 0)  to (1, .5);   
    \draw (.5, 0)  to (1, -0.5); 

    \draw (-1, .5)  to (1, .5);   
    \draw (1, .5)  to (1, -0.5); 
    \draw (1,-0.5) to (-1,-0.5);  
    \draw (-1,-0.5) to (-1, .5);

    \fill (-.5,0) circle(.1);
    \fill (.5,0) circle(.1);
    \fill (-1,0.5) circle(.1);
    \fill (-1,-0.5) circle(.1);
    \fill (1,0.5) circle(.1);
    \fill (1,-0.5) circle(.1);

    \node[align=center] (H4) at (0,-1) {$H_{4}$ \\ 4 cycles};
\end{scope}

\begin{scope}[shift={(10,2)}]
    \draw (0,0) to (-1,.5);
    \draw (0,0) to (-1,-0.5);
    \draw (0,0) to (1,.5);
    \draw (0,0) to (1,-0.5);
    
    \draw (0,2) to (-1,0.5);
    \draw (0,2) to (-1,-0.5);
    \draw (0,2) to (1,.5);
    \draw (0,2) to (1,-0.5);
    \draw (0,2) to (2,0);

    \draw (-1, .5) to (1, .5);
    \draw (-1, .5) to (-1, -0.5);

    \draw (1, .5) to (2, 0);

    \draw (2, 0) to (1, -0.5);
    \draw (1, -0.5) to (-1, -0.5);

    \fill (0,0) circle(.1);
    \fill (1,.5) circle(.1);
    \fill (0,2) circle(.1);
    \fill (-1,.5) circle(.1);
    \fill (-1,-0.5) circle(.1);
    \fill (2,0) circle(.1);
    \fill (1,-0.5) circle(.1);

    \node[align=center] (G5) at (0,-1) {$G = G_5$\\8 cycles};
\end{scope}

\begin{scope}[shift={(10,-2)}]
    \draw (-.5,0)   to (-1,.5);    
    \draw (-.5,0)   to (-1,-0.5);  
    \draw (-.5,0)   to (.5,0);       

    \draw (.5, 0)  to (1, .5);   
    \draw (.5, 0)  to (1, -0.5); 

    \draw (-1, .5)  to (1, .5);   
    \draw (1, .5)  to (1, -0.5); 
    \draw (1,-0.5) to (-1,-0.5);  
    \draw (-1,-0.5) to (-1, .5);

    \fill (-.5,0) circle(.1);
    \fill (.5,0) circle(.1);
    \fill (-1,0.5) circle(.1);
    \fill (-1,-0.5) circle(.1);
    \fill (1,0.5) circle(.1);
    \fill (1,-0.5) circle(.1);
    
    \draw (0,2) to (-1,0.5);
    \draw (0,2) to (-1,-0.5);
    \draw (0,2) to (-.5,0);
    \draw (0,2) to (1,-0.5);
    \draw (0,2) to (1,0.5);

    \fill (0,2) circle(.1);

    \node[align=center] (H5) at (0,-1) {$H = H_5$ \\ 8 cycles};
\end{scope}
\draw[->,shorten >=30pt, shorten <=10pt] (G3.east) -- (G4.west);
\draw[->,shorten >=20pt, shorten <=20pt] (G4.east) -- (G5.west);
\draw[->,shorten >=40pt, shorten <=10pt] (G3.east) -- (H4.west);
\draw[->,shorten >=20pt, shorten <=20pt] (H4.east) -- (H5.west);
\end{tikzpicture}
    \caption{Non-isomorphic codegree pair of graphs $G$ and $H$. For $k=3, 4, 5$, the graphs $G_k$ and $H_k$ denote the subgraphs obtained from $G$ and $H$, respectively, by degree filtration at level $k$. They coincide at $k=3$, but evolve differently for higher levels. The full graphs $G$ and $H$ share the same second smallest eigenvalue. Note that $H_1(G_4) = H_1(H_4) = 4$, and $H_1(G_5)=H_1(H_5)=8$.}
    \label{fig:eigenval_dist2}
\end{figure}

\begin{figure}[ht] 
    \centering
    \begin{subfigure}{0.48\textwidth} 
        \centering
        \includegraphics[width=\textwidth]{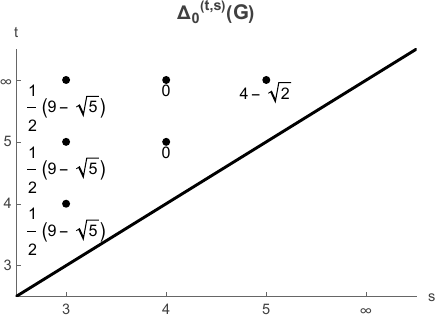}
    \end{subfigure}
    \hfill 
    \begin{subfigure}{0.48\textwidth} 
        \centering
        \includegraphics[width=\textwidth]{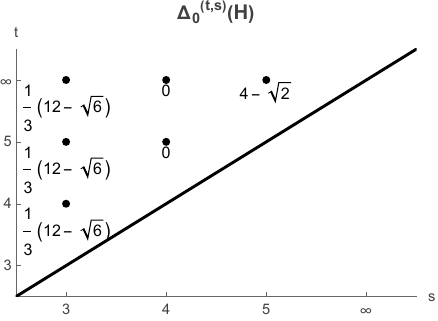}
    \end{subfigure}
    \caption{Persistent Laplacian diagrams of $G$ and $H$ in Figure \ref{fig:eigenval_dist2} with signature being their second smallest eigenvalues}
    \label{fig:pld}
\end{figure}

In Figure~\ref{fig:eigenval_dist2}, the graphs $G$ and $H$ have the same degree sequence $(3, 4, 4, 4, 4, 4, 5)$, and the second smallest eigenvalue $4 - \sqrt{2}$. In degree filtration, both graphs have the same number of simplicial homology cycles at each degree level and also the same persistence diagrams such that $B_0 = \{(3, \infty)\}$ and $B_1 = \{(4, \infty), (5, \infty)\}$ with multiplicity $4$ for each element in $B_1$. They have different persistent Laplacians between degree 3 and 4, and their second smallest eigenvalues $s_{gap}$ are $s_{gap}\left(\Delta_0^{4,3}(G)\right) = \frac{9 - \sqrt{5}}{2}$, and $s_{gap}\left(\Delta_0^{4,3}(H)\right) = \frac{12-\sqrt6}{3}$, see Figure \ref{fig:pld}.

\begin{remark}
    For $a<b<c$, we expect a connection between $\Delta_0^{a,b}$ and $\Delta_0^{a,c}$ in a simple circumstance, for example when all deaths occur at $\infty$.
\end{remark}

\subsection{Cospectral Graphs} \label{sec: cospectral graphs}

There exist pairs of graphs that are both codegree and cospectral. A well-known example is the Shrikhande graph and the \(4 \times 4\) rook's graph. Both graphs have \(16\) vertices and \(48\) edges, with each vertex having degree \(6\).

To define these graphs, label the vertices by ordered pairs \((i,j)\) where \(i,j \in \mathbb{Z}/4\mathbb{Z}\).  
In the Shrikhande graph $G_S$, two vertices \((i,j)\) and \((i',j')\) are adjacent if and only if
\[
(i-i', j-j') \in \{ (\pm 1, 0), (0, \pm 1), (\pm 1, \pm 1) \}.
\]
For the \(4 \times 4\) rook's graph $G_R$, two distinct vertices \((i,j)\) and \((i',j')\) are adjacent if and only if \(i = i'\) or \(j = j'\).

\begin{figure}[h]
    \centering
    \begin{minipage}{0.48\textwidth}
        \centering
\begin{tikzpicture}[ scale=1.5, every node/.style={circle, draw, fill=white, inner sep=0.8pt},
smaller nodes nodestyle/.style={font=\scriptsize},baseline={(current bounding box.center)} ] 
\foreach \i in {0,...,3}{
\foreach \j in {0,...,3}{
\node (v\i\j) at (\i,\j) {}; } } 
\foreach \i in {0,...,3}{ \foreach \j in {0,...,3}{ \foreach \di/\dj in {1/0, 3/0, 0/1, 0/3, 1/1, 3/3}{ \pgfmathtruncatemacro{\ii}{mod(\i+\di,4)} \pgfmathtruncatemacro{\jj}{mod(\j+\dj,4)} \pgfmathtruncatemacro{\fromindex}{4*\i+\j} \pgfmathtruncatemacro{\toindex}{4*\ii+\jj} \ifnum\fromindex<\toindex 
\ifnum\di=3 
\draw[bend left=10] (v\i\j) to (v\ii\jj); \else\ifnum\dj=3
\draw[bend left=20] (v\i\j) to (v\ii\jj); \else 
\draw (v\i\j) -- (v\ii\jj); \fi\fi \fi } } } 
\foreach \i in {0,...,3}{ \foreach \j in {0,...,3}{ \fill (\i,\j) circle(.05); } } 
\end{tikzpicture}
\end{minipage}%
\hspace{0.02\textwidth}%
\begin{minipage}{0.48\textwidth}
        \centering
\begin{tikzpicture}[scale=1.5, every node/.style={circle, draw, fill=white, inner sep=0.8pt},
nodestyle/.style={font=\scriptsize}, baseline={(current bounding box.center)} ] 
\foreach \i in {0,...,3}{ \foreach \j in {0,...,3}{ \node (v\i\j) at (\i,\j) {}; } }
\foreach \i in {0,...,3}{ \foreach \j in {0,...,3}{ \foreach \di/\dj in {1/0,2/0,3/0,0/1,0/2,0/3}{ \pgfmathtruncatemacro{\ii}{mod(\i+\di,4)} \pgfmathtruncatemacro{\jj}{mod(\j+\dj,4)} \pgfmathtruncatemacro{\fromindex}{4*\i+\j} \pgfmathtruncatemacro{\toindex}{4*\ii+\jj} \ifnum\fromindex<\toindex
\ifnum\di=2 \draw[bend left=20] (v\i\j) to (v\ii\jj); \else\ifnum\dj=2 \draw[bend left=20] (v\i\j) to (v\ii\jj); \else\ifnum\di=3 \draw[bend left=20] (v\i\j) to (v\ii\jj); \else\ifnum\dj=3 \draw[bend left=20] (v\i\j) to (v\ii\jj); \else 
\draw (v\i\j) -- (v\ii\jj); \fi\fi\fi\fi \fi } } }
\foreach \i in {0,...,3}{ \foreach \j in {0,...,3}{ \fill (\i,\j) circle(.05); } } 
\end{tikzpicture}
\end{minipage}
    \caption{Shrikhande graph $G_S$ (left) and $4\times 4$ rook's graph $G_R$(right).}
    \label{fig:ShrikhandeRook}
\end{figure}

These graphs are both codegree and cospectral, meaning that they share the same degree sequence and have identical eigenvalues for their graph Laplacians. However, their eigenvectors differ. Therefore, one may consider the geometry-based signatures $s_{\textrm{geo},2}[v]$ associated with $v = (1,0,\ldots, 0)$, since
$s_{\textrm{geo},2}\left(\Delta^{\infty, 6}_0(G_S)\right) \simeq 7.3$ and
$s_{\textrm{geo},2}\left(\Delta^{\infty, 6}_0(G_R)\right) \simeq 10.6$
are distinct.

\section{Examples of Signatures}\label{subsec:sign}

We now introduce examples of additional variables, referred to as \emph{signatures} of the persistent Laplacian. Examples of such signatures include classical eigenvalue-based quantities such as the second smallest eigenvalue and spectral entropy, and to overcome their limitations under cospectrality, we additionally introduce an eigenvector-based signature. Let $N := n_q^{K_d}-n_q^{K_b}$ be the size of the $q$-th persistent Laplacian $\Delta_q^{K_d,K_b}$. All suggested signatures require the computation of eigenpairs.

\subsection{Eigenvalues-based signatures}

A natural class of signatures arises from the spectrum of the persistent Laplacian. For a simplicial pair $(K_b \hookrightarrow K_d)$ with $b < d$, recall that the $q$-th persistent Laplacian $\Delta^{K_d,K_b}_q$ is a real symmetric positive semi-definite matrix. Hence, all its eigenvalues are real and nonnegative. Let
\[
0 = \lambda_1 \le \lambda_2 \le \cdots \le \lambda_{N}
\]
denote the eigenvalues of $\Delta^{K_d,K_b}_q$ in non-decreasing order.

\begin{remark}
    Spectral information from the entire spectrum was incorporated into vector representations via truncation, suggesting that eigenvalues capture important structural features \cite{davies2023persistent}.
\end{remark}

\paragraph{Second smallest eigenvalue}

The second smallest eigenvalue $\lambda > 0$ plays a role analogous to the spectral gap in classical graph Laplacian theory~\cite{Chung97, fiedler1973algebraic, spielman2012spectral}. This interpretation extends naturally to persistence diagrams and persistent Laplacians. In~\cite{chen2022persistent}, the authors employed various spectral quantities, including the second smallest eigenvalue, to analyze biomolecular data. Broadly speaking, $\lambda$ measures the \emph{strength of connectivity}, or the degree of coupling, among the $q$-chains that persist over the interval $[b,d]$. The map
\[
s_{\mathrm{gap}}(\Delta^{K_d,K_b}_q) := \lambda
\]
defines an admissible signature in the sense of Definition \ref{def:signatures}. $\lambda$ is bounded above by $2h^{K_b}$ where $h^{K_b}$ is the Cheeger constant \cite{memoli2022persistent}.

\paragraph{Spectral entropy.}
Another eigenvalue-based signature is the \emph{spectral entropy}, which measures the spread or disorder of the spectrum \cite{de2016spectral, nicolini2018thermodynamics, sun2021eigenvalue}. Define the normalized eigenvalue weights
\[
p_i := \frac{\lambda_i}{\sum_{j=1}^{N} \lambda_j}, \qquad \sum_i p_i = 1,
\]
with the convention that $p_i = \frac{1}{N}$ when all eigenvalues vanish.
The spectral entropy of the persistent Laplacian is
\[
s_{\mathrm{ent}}(\Delta^{K_d,K_b}_q)
  := - \sum_{i=1}^{N} p_i \log p_i.
\]
$s_{\mathrm{ent}}$ is bounded by $\log N$, which makes $s_{\mathrm{ent}}$ another admissible signature for use in the persistent Laplacian diagram and image.

\paragraph{Trace and spectral moments.}
The trace of the persistent Laplacian provides another simple eigenvalue-based signature. Since $\Delta^{K_d,K_b}_q$ is symmetric and positive semi-definite, its trace is equal to the sum of its eigenvalues:
\[
s_{\mathrm{tr}}(\Delta^{K_d,K_b}_q)
    := \mathrm{tr}(\Delta^{K_d,K_b}_q)
    = \sum_{i=1}^{N} \lambda_i.
\]
Thus $s_{\mathrm{tr}}$ measures the total spectral mass of the persistent Laplacian over the interval $[b,d]$. While spectral entropy records how this mass is distributed across the spectrum, the trace records its overall scale. More generally, one may use spectral moment signatures~\cite{estrada1996spectral, godsil2013algebraic, mowshowitz2014spectral}
\[
s_{\mathrm{mom},r}(\Delta^{K_d,K_b}_q)
    := \left(\sum_{i=1}^{N} \lambda_i^r\right)^{1/r},
    \qquad r \ge 1,
\]
These signatures interpolate between different summaries of the eigenvalue distribution. The case $r=1$ recovers the trace, while larger values of $r$ emphasize the contribution of large eigenvalues. On any fixed finite data class, or after applying a normalization that uniformly bounds the spectrum, these quantities are admissible signatures in the sense of Definition \ref{def:signatures}.

\subsection{Geometry-based signatures}\label{sec:l1l2}
While eigenvalue-based signatures capture global spectral quantities of the persistent Laplacian, they fail to distinguish matrices that are \emph{cospectral}, meaning that two different Laplacians have the same eigenvalues. Two persistent Laplacians may share identical eigenvalues but differ significantly in the 
geometric structure encoded by their eigenvectors. To overcome this limitation, we introduce a novel geometry-based signature that incorporates eigenvector information.

Cospectrality is a classical phenomenon in spectral graph theory, and the same issue appears for persistent Laplacians. If two operators $\Delta^{K_d,K_b}_q$ have the same spectra, all eigenvalue-based signatures necessarily agree. Therefore, such signatures fail to separate different geometric or combinatorial behaviors present in the persistent pair $(K_b \hookrightarrow K_d)$. This motivates the use of eigenvector-based quantities.

\paragraph*{Geometric eigenvector profile signature}

Let $\left\{ \lambda_i \right\}_{i=1}^N$ be the spectrum of $\Delta^{K_d,K_b}_q$, $E_{\lambda_i}$ be the eigenspace corresponding to $\lambda_i$. For any vector $v\in\R^N$ and any subspace $E \subset \R^N$, let $\textrm{Proj}_{E}(v)$ denote the orthogonal projection $v$ onto $E$. Fix a reference vector $v \in \R^{N}$ and a positive integer $p$. We define the \emph{geometric eigenvector profile signature} $s_{\mathrm{geo}}[v]$ associated with $v$ by
\[
s_{\mathrm{geo},p}[v](\Delta^{K_d,K_b}_q) := \sum_{E_{\lambda_i}} \|\textrm{Proj}_{E_{\lambda_i}}(v)\|_p.
\]
Since $\|\mathrm{Proj}_{E_{\lambda_i}}(v)\|_p$ is independent of the choice of basis for each $E_{\lambda_i}$, the quantity $s_{\mathrm{geo}}[v]$ is well-defined. Since $s_{\mathrm{geo},p}[v]$ is bounded by $\|v\|_1 < \infty$, $s_{geo}$ is admissible. When the vector $v$ is clear from context, we omit it from the notation. Note that if $v$ happens to be a kernel of Laplacian, then $s_{\mathrm{geo}}=0$. Also, we would like to avoid $v$ being an eigenvector. This provides a mechanism to distinguish cospectral operators even when the eigenvalues agree, the geometry of the eigenvectors often differs substantially, and the quantity $s_{\mathrm{geo}}$ reflects this.


\section{Exploratory Results Across Signatures}
\label{sec:additional-expeiment}

Before conducting the repeated experiments reported in Tables~\ref{tab:mnist_accuracy_metrics} and~\ref{tab:persistent_laplacian_metrics}, we performed an exploratory comparison over a broader collection of persistent-Laplacian signatures. This comparison included the smallest positive eigenvalue, spectral entropy, trace, and geometric \(L_1\) signatures. It was used to identify promising signatures and to examine the effects of filtration and image resolution.

\begin{table}[t]
\centering
\caption{Exploratory MNIST accuracy point estimates across PLI signatures. The PL column denotes the reproduced persistent-Laplacian baseline and does not depend on image resolution.}
\label{tab:old_mnist}
\small
\setlength{\tabcolsep}{4pt}
\begin{tabular}{@{}lccccc@{}}
\toprule
\(\boldsymbol{(r_f,r_i)}\) &
\textbf{PL} &
\(\boldsymbol{\lambda_{>0,\min}}\) &
\textbf{Spectral entropy} &
\textbf{Trace} &
\textbf{Geometric \(L_1\)} \\
\midrule
\((4,10)\)  & \multirow{3}{*}{0.6566} & 0.6255 & 0.5607 & 0.5972 & 0.3986 \\
\((4,20)\)  &                         & 0.6189 & 0.5589 & 0.5965 & 0.4014 \\
\((4,30)\)  &                         & 0.6179 & 0.5578 & 0.5971 & 0.3982 \\
\midrule
\((8,10)\)  & \multirow{3}{*}{0.8199} & 0.8055 & 0.8228 & 0.8345 & 0.5775 \\
\((8,20)\)  &                         & 0.8021 & 0.8268 & 0.8436 & 0.5732 \\
\((8,30)\)  &                         & 0.8024 & 0.8252 & 0.8397 & 0.5735 \\
\midrule
\((12,10)\) & \multirow{3}{*}{0.8326} & 0.8400 & 0.8581 & 0.8781 & 0.6508 \\
\((12,20)\) &                         & 0.8379 & 0.8625 & \textbf{0.8850} & 0.6558 \\
\((12,30)\) &                         & 0.8424 & 0.8611 & 0.8837 & 0.6504 \\
\bottomrule
\end{tabular}
\end{table}

\begin{table}[t]
\centering
\caption{Exploratory QM7 MAE point estimates across PLI signatures. The PL column denotes the reproduced persistent-Laplacian baseline and does not depend on image resolution.}
\label{tab:old_qm7}
\small
\setlength{\tabcolsep}{4pt}
\begin{tabular}{@{}lccccc@{}}
\toprule
\(\boldsymbol{(r_f,r_i)}\) &
\textbf{PL} &
\(\boldsymbol{\lambda_{>0,\min}}\) &
\textbf{Spectral entropy} &
\textbf{Trace} &
\textbf{Geometric \(L_1\)} \\
\midrule
\((4,10)\)  & \multirow{3}{*}{60.5902} & 65.1954 & 41.2708 & 40.0596 & 65.8822 \\
\((4,20)\)  &                          & 65.1630 & 41.2602 & 39.9395 & 66.6892 \\
\((4,30)\)  &                          & 66.0509 & 41.4091 & 40.0651 & 64.5358 \\
\midrule
\((8,10)\)  & \multirow{3}{*}{54.7390} & 56.9239 & 39.5125 & 38.2591 & 56.7884 \\
\((8,20)\)  &                          & 56.8367 & 39.4651 & 38.0948 & 56.2647 \\
\((8,30)\)  &                          & 57.5178 & 39.5530 & 38.4441 & 57.9197 \\
\midrule
\((12,10)\) & \multirow{3}{*}{50.1028} & 52.2036 & 40.8944 & 38.2505 & 52.0541 \\
\((12,20)\) &                          & 52.9954 & 42.0473 & 37.9288 & 52.2678 \\
\((12,30)\) &                          & 52.0839 & 40.8764 & \textbf{37.4249} & 52.6692 \\
\bottomrule
\end{tabular}
\end{table}

For MNIST, trace gives the strongest PLI result at filtration resolutions \(r_f=8\) and \(r_f=12\), reaching an accuracy of \(0.8850\) at \((r_f,r_i)=(12,20)\). Spectral entropy is generally the next strongest signature at these resolutions. For QM7, trace attains the lowest MAE in every tested resolution configuration, with its best point estimate of \(37.4249\) at \((r_f,r_i)=(12,30)\). Based on these exploratory results, the controlled repeated-run experiments in the main text evaluate trace and spectral entropy.

These results also suggest that filtration resolution has a larger effect than image resolution in the tested configurations. Increasing \(r_f\) generally improves the trace-based representation, whereas changing \(r_i\) while holding \(r_f\) fixed produces comparatively smaller differences. This is consistent with the roles of the two resolutions: filtration resolution controls how finely the evolution of the underlying complex is sampled, whereas image resolution controls only the discretization of the resulting persistent-Laplacian surface. These observations are exploratory and are not intended as formal conclusions about optimal resolution selection.

\section{Algorithms for the Persistent Laplacian}
\label{app:pl-algorithms}

We summarize the computational steps used to construct a persistent Laplacian along an inclusion $K\hookrightarrow L$. The presentation follows the Schur-complement, or Kron-reduction, viewpoint of~\cite{memoli2022persistent}.

\subsection{Schur complement and Kron reduction}

Recall that the up part of the $q$-th persistent Laplacian is
\[
\Delta^{L,K}_{q,\mathrm{up}}
  = \partial^{L,K}_{q+1}\circ
    \bigl(\partial^{L,K}_{q+1}\bigr)^{*},
\]
and that the full operator is
\[
\Delta^{L,K}_{q}
  = \Delta^{L,K}_{q,\mathrm{up}}
    + \bigl(\partial_q^K\bigr)^{*}\circ\partial_q^K.
\]
Let $n_q^L=|S_q^L|$ and $n_q^K=|S_q^K|$, and write
$\Delta^L_{q,\mathrm{up}}=\partial^L_{q+1}\circ
(\partial^L_{q+1})^*$ for the usual up-Laplacian on $L$.
In the unweighted case, or after passing to a weighted orthonormal basis, order the
$q$-simplices so that those in $K$ appear first. Then
\[
\Delta^L_{q,\mathrm{up}}
=
\begin{bmatrix}
A & B\\
B^{\mathsf T} & D
\end{bmatrix},
\qquad
A\in\R^{n_q^K\times n_q^K},
\]
where $D$ is indexed by the $q$-simplices in $L\setminus K$.
Theorem~4.6 of~\cite{memoli2022persistent} identifies the up-persistent
Laplacian with the generalized Schur complement
\[
\Delta^{L,K}_{q,\mathrm{up}}
  = \Delta^L_{q,\mathrm{up}}/D
  = A-BD^{\dagger}B^{\mathsf T},
\]
where $D^{\dagger}$ is the Moore--Penrose pseudoinverse. This is precisely a
Kron reduction, familiar from electrical-network theory, lifted to the
persistent setting~\cite{dorfler2012kron,kron1939tensor}.

\subsection{Algorithmic construction}

Let $B^L_{q+1}$ and $B^K_q$ be the boundary matrices of
$\partial^L_{q+1}$ and $\partial^K_q$, respectively, in the canonical bases,
and let $W^L_{q+1}$ and $W^K_q$ be the corresponding diagonal weight matrices.
For a matrix $M\in\R^{m\times n}$ and index sets $I\subseteq[m]$ and
$J\subseteq[n]$, let $M(I,J)$ denote the submatrix with rows in $I$ and
columns in $J$; the notation $M(:,J)$ and $M(I,:)$ has the usual meaning.

Using the same $q$-simplex ordering and following Lemma~3.4
of~\cite{memoli2022persistent}, define
\[
D^L_{q+1}
  :=B^L_{q+1}\bigl([n_q^L]\setminus[n_q^K],:\bigr).
\]
Column reduction gives $R^L_{q+1}=D^L_{q+1}Y$ for an invertible matrix $Y$.
Let $I$ index the zero columns of $R^L_{q+1}$. The associated columns of $Y$
form a basis for
\[
C^{L,K}_{q+1}
  =\{c\in C^L_{q+1}:\partial^L_{q+1}c\in C^K_q\}.
\]
Set
\[
Z:=Y(:,I),
\qquad
B^{L,K}_{q+1}:=
  (B^L_{q+1}Y)([n_q^K],I),
\]
so that $B^{L,K}_{q+1}$ represents $\partial^{L,K}_{q+1}$ in this basis.
If $I=\emptyset$, then $C^{L,K}_{q+1}=\{0\}$ and the up term vanishes.
Otherwise,
\[
\Delta^{L,K}_{q,\mathrm{up}}
=B^{L,K}_{q+1}
 \bigl(Z^{\mathsf T}(W^L_{q+1})^{-1}Z\bigr)^{-1}
 (B^{L,K}_{q+1})^{\mathsf T}(W^K_q)^{-1}.
\]
Adding the down-Laplacian gives
\[
\Delta^{L,K}_q
=\Delta^{L,K}_{q,\mathrm{up}}
 +W^K_q(B^K_q)^{\mathsf T}(W^K_{q-1})^{-1}B^K_q.
\]
This expression is invariant under the choice of basis for
$C^{L,K}_{q+1}$.

\subsection{Time and memory complexity}

We state dense-arithmetic upper bounds; sparsity can substantially reduce the
cost. Reducing the $(n_q^L-n_q^K)\times n_{q+1}^L$ matrix $D^L_{q+1}$ costs
\[
O\bigl((n_q^L-n_q^K)(n_{q+1}^L)^2\bigr).
\]
Forming $B^L_{q+1}Y$ costs $O(n_q^L(n_{q+1}^L)^2)$, while forming and
inverting the restricted Gram matrix costs at most $O((n_{q+1}^L)^3)$.
For fixed $q$, the down-Laplacian can be assembled in
$O((n_q^K)^2)$ time by exploiting the sparsity of simplicial boundary
matrices. Thus, a dense upper bound for constructing one
$\Delta_q^{L,K}$ is
\begin{equation}
\label{eq:time-complexity}
O\Bigl(
n_q^L(n_{q+1}^L)^2
+(n_{q+1}^L)^3
+(n_q^K)^2
\Bigr).
\end{equation}
The corresponding dense working memory is
\begin{equation}
\label{eq:memory-complexity}
O\Bigl(
n_q^L n_{q+1}^L
+(n_{q+1}^L)^2
+(n_q^K)^2
\Bigr),
\end{equation}
accounting for the boundary data, change-of-basis matrix, and output
persistent Laplacian. Fast matrix multiplication, sparse linear algebra, and
reuse of Schur complements along a filtration can improve these bounds; see
Section~5 of~\cite{memoli2022persistent}.

For the complete PLI pipeline, let $P$ be the number of selected birth--death
pairs and let $m_q$ be the largest persistent-Laplacian dimension among them.
A full dense eigendecomposition costs $O(m_q^3)$ time and $O(m_q^2)$ memory per
pair, hence $O(Pm_q^3)$ spectral time if all $P$ pairs are processed. The
matrices may be streamed, so the spectral working memory remains $O(m_q^2)$;
retaining all of them would instead require $O(Pm_q^2)$ storage. On a complex
with $n$ vertices, $m_q\leq\binom{n}{q+1}$. In the graph-based regime used in
our experiments, $m_0\leq n$, while $m_1$ is at most the number of edges.
Consequently, a sparse graph with $m_1=O(n)$ has $O(n^2)$ dense spectral
memory and $O(n^3)$ full-eigendecomposition time per pair, whereas the dense
$q=1$ worst case reaches $O(n^4)$ memory and $O(n^6)$ time.

The required spectral work depends on the signature. The trace is available
from the diagonal in $O(m_q)$ time once the persistent Laplacian is formed and
does not require eigendecomposition. The smallest positive eigenvalue can be
estimated with a sparse partial eigensolver, with cost governed by matrix
sparsity and convergence. Spectral entropy and the geometric eigenvector
profile generally require a larger portion of the spectrum. Thus, in the
present implementation, persistent-Laplacian construction and spectral
analysis dominate Gaussian smoothing and pixel integration; the latter do not
remove the scaling bottleneck above.

\section{Resources and additional experimental detail}\label{app:resource}

The computational experiments were performed on a high-performance workstation equipped with 16 Intel Xeon Gold 6426Y CPU, 256GB of RAM, and an NVIDIA L40S GPU. The running time is described in Table~\ref{tab:computation_runtimes}. These measurements complement the asymptotic construction and spectral costs in Appendix~\ref{app:pl-algorithms}.

\begin{table}[t]
\centering
\caption{Average Computational Running Times (in seconds) for MNIST and QM7 Datasets}
\label{tab:computation_runtimes}
\begin{tabular}{cccc}
\hline
\textbf{Dataset} & \textbf{Filtration Resolution} & \textbf{Image Resolution} & \textbf{Average Runtime (s)} \\ \hline
\multirow{7}{*}{\textbf{MNIST}} & 4 & All (10, 20, 30) & 261 \\ \cline{2-4} 
                                & \multirow{3}{*}{8} & 10 & 960 \\
                                &                    & 20 & 989 \\
                                &                    & 30 & 1038 \\ \cline{2-4} 
                                & \multirow{3}{*}{12} & 10 & 2088 \\
                                &                    & 20 & 2114 \\
                                &                    & 30 & 2162 \\ \hline
\multirow{3}{*}{\textbf{QM7}}   & 4                  & All (10, 20, 30) & 12 -- 24 \\
                                & 8                  & All (10, 20, 30) & 48 -- 54 \\
                                & 12                 & All (10, 20, 30) & $\sim$ 114 \\ \hline
\end{tabular}
\end{table}

For downstream prediction, we used standard multilayer perceptrons implemented in PyTorch. For MNIST, the dataset was split into training, validation, and test sets with ratio $0.64:0.16:0.20$. The classifier consisted of three hidden layers of width $128$, each followed by a ReLU activation. For QM7, the split ratio is $0.72:0.08:0.20$. The regressor was also implemented as an MLP, with four hidden layers of width $128$, each followed by a ReLU activation.

For both tasks, optimization was performed using Adam with learning rate $10^{-3}$ and weight decay $10^{-5}$. We used a cosine annealing learning-rate scheduler with minimum learning rate $10^{-5}$ over $200$ epochs.


\end{document}